\documentclass[12pt]{article}
     \usepackage{latexsym}
     \usepackage{amsmath}
     \usepackage{amssymb}

\begin{document}
\begin{center}
{\huge\bf On maps which preserve equality of distance in
$F^{*}$-spaces
 \footnote{This paper is supported by The
National Natural Science Foundation of China (10571090) and The
Research Fund for the Doctoral Program of Higher Education
(20010055013).}} \baselineskip 5mm \vskip 5mm  Dongni ,Tan

{\small School of Mathematics Science, Nankai University,
Tianjin300071, China}

{\small  E-mail: 0110127@mail.nankai.edu.cn} \baselineskip 5mm
\vskip 5mm
\begin{minipage}{5in}
\baselineskip 5mm {\bf Abstract}\quad In order to generalize the
results of Mazur-Ulam and Vogt, we shall prove that any map $T$
which preserves equality of distance with $T(0)=0$ between two
$F^{*}$-spaces without surjective condition is linear. Then , as a
special case linear isometries are characterized through a simple
property of their
range.\\

\end{minipage}
\end{center}
\baselineskip 7.5mm
\section{ Preliminaries and introduction}
\indent \indent Recall from [1] that a non-negative function $\|.\|$
defined on a linear space $E$ is called an $F$-norm provided
\par(i) $\|x\|=0$ iff $x=0$;\par(ii) $\|ax\|=\|x\|$ for all $a$ with
$|a|=1$;\par (iii)  $ \|x+y\|\leqslant \|x\|+\|y\|$;\par(iv)
$\|a_nx\|\rightarrow 0 $ provided $a_n\rightarrow 0$; \par(v)
$\|ax_n\|\rightarrow 0 $ provided $x_n \rightarrow 0$.\par An
$F$-norm $\|.\|$ induces a transitive invariant distance $d$ by
$$d(x,y)=\|x-y\| \quad \forall\,x,y\in E.$$
A linear space $E$ with an $F$-norm $\|.\|$ is said to be an
$F^*$-space. The complete $F^*$-space is called $F$-space.\\\indent
Let $E$ =$(E,\,\|.\|)$, $F=(F,\, \|.\|)$ be two $F^*$-spaces, Let
$R_0^+$ denote the set of non-negative real numbers. We say that a
map $T:E\rightarrow F$ preserves equality of distance iff there
exists a function $\phi:R_0^+ \rightarrow R_0^+$ such that $$
\|Tx-Ty\|=\phi\|x-y\|,\quad  \forall\,x,y\in E.$$ The function
$\phi$ is called the gauge function for $T$, and this definition
equivalently means that $$u,\,v,\,z,\,w\, \in E,\,\,
\|u-v\|=\|z-w\|\Rightarrow \|Tu-Tv\|=\|Tz-Tw\|.$$In particular, when
$\phi$ is the identity function, mapping $T$ is an isometry. Such
mappings were studied by Schoenberg [2] and by Von Neumann and
Schoenberg [3] for Hilbert spaces.\\\indent The classical theorem of
Mazur-Ulam [4] states that an onto isometry between two real normed
spaces is affine, Charzy$\acute{n}$ski [5] and Rolewicz [6] have
shown, respectively, that surjective isometries of
finite-dimensional $F$-spaces and of locally bounded spaces with
concave norm are also linear. Ding and huang [8] showed that
Rolewicz $^,s$ result also holds in locally midpoint constricted
$F$-spaces. More generally, Vogt [7] has shown that equality of
distance preserving maps between two real normed spaces are linear.
The hypothesis of surjectivity of the maps was required in the
results mentioned above.
\\\indent The present paper will extend the result of Vogt to a large
class of $F^*$-spaces including all $p$-normed spaces$(0<p\leq1)$
 without surjective
condition, but needing a simple peculiar property of the maps$^,$
range. The proof of the main result here depends on the technique of
Vogt and all the spaces mentioned in this paper are assumed to be
real.
\section{ Main results}
\qquad The following lemma is in metric space theory,  which was
adapted by Vogt in [7] and similar to one stated by Aronszajn[9].\\
\textbf{Lemma 2.1.} \,Let $(M,d)$ be a bounded metric space. Suppose
there exists an element $m$ in $M$, a surjective isometry $V:
M\rightarrow M $, a constant $K>1$, such that for all $x$ in $M$,
$d(Vx,x)\geq Kd(m,x)$, Then $m$ is a fixed point for every
surjective isometry $S: M\rightarrow M.$\\
\textbf{Proof.} \,Since metric space isometries are injective,
$V^{-1}\, \mbox{and}\, S^{-1}$ exist, and $V,\, S,\, V^{-1}$ and
$S^{-1}$ are bijective isometries of $M$ together with arbitrary
(finite) compositions of them.\\\indent Define a sequence of
isometries $V_n : M\rightarrow M$ and elements $m_n$ in $M$  indexed
by the integers $n\geq 1$,\par Let $$ V_1=V,  \qquad \qquad \qquad
m_1=m,$$
$$ V_2=SVS^{-1}, \qquad \qquad  m_2=Sm,$$ $$ \qquad \qquad
V_{n+1}=V_{n-1}V_n(V_{n-1})^{-1},\qquad m_{n+1}=V_{n-1}m_n,\quad
n\geq 2.$$\indent Each $V_n$ is a bijective, invertible isometry of
$M$ and a straight-forward induction yields:$$ d(V_nx,\,x)\geq K
d(m_n,\,x)\quad x\in M,\,n\geq 1.\eqno{(1)}$$ If we let $x=m_{n+1}$
in (1), we obtain $$
d(m_{n+2},\,m_{n+1})=d(V_nm_{n+1},\,m_{n+1})\geq
Kd(m_n,\,m_{n+1})=Kd(m_{n+1}\, m_n).$$ Another induction gives
$$ d(m_{n+2},\,m_{n+1})\geq K^nd(m_2,m_1) \qquad \forall n\geq 1.$$
\indent Since $M$ is a bounded metric space, there exists a positive
number $N$ such that $$  d(m_{n+2},\,m_{n+1})\leq N \qquad \forall
n\geq 1.$$ Hence
$$ d(m_2,\,m_1)\leq N/K^n \qquad \forall n\geq 1.$$ Since $K>1$,we have
$d(m_2,\,m_1)=0.$ This implies $$ Sm=m_2=m_1=m.\qquad \Box$$

\par
Rassias [12] shows that the ratio $\|2x\|/\|x\|$ plays an important
role in the generalizations of the Mazur-Ulam theorem and we shall
show that the ratio $\|2x\|/\|x\|$ is also of great importance in
the following statement.
 \par \textbf{Theorem 2.2.} \,Let $E,F$ be $F^*$-spaces, and there exists a positive number $r$
such that $ \alpha_F(r)=\inf\{\frac{\|2x\|}{\|x\|}:\, x\in F,\,
0<\|2x\|\leq r\}>1$, Let $T: E\rightarrow F \,\mbox {with}\, T(0)=0$
be a continuous map which preserves equality of distance, i.e. there
exists $\phi:R_0^+ \rightarrow R_0^+$ such that $$
\|Tx-Ty\|=\phi\|x-y\|,\quad  \forall\,x,y\in E.$$ and satisfy $$
2Tx-Ty\in T(E),\quad \forall \, x, y \in E.$$ Then T is linear.
\par \textbf{Proof.}\, Since $T$ is continuous, there is a $\delta>0$
such that for any $a, b\in E$ satisfying $\|\frac {a-b}{2}\|\leq
\delta$, we have $\|T (\frac {a-b}{2})\|\leq \frac {r}{4}.$ \par
\textbf{Step 1.}\,For any $a,b \in E$ satisfying $\|\frac
{a-b}{2}\|\leq\delta.$ \par Let $m=T(\frac{a-b}{2})$ and define $M
:= \{ y\in T(E): \|y\|=\|2m-y\|\leq2\|m\|\leq \frac {r}{2}\}$ and
$V: M\rightarrow M \, \mbox {by}\, V(y)=2m-y.$ Then we get \par
(i)\, $M\neq\varnothing \,\mbox {since}\, m\in M.$\par  (ii)\,$M$ is
a bounded metric space.
\par (iii)\,$V$ is an isometry from $M$ onto $M$since $
V^2=id_M$ and $\|V(y_1)-V(y_2)\|=\|2m-y_1-2m+y_2\|=\|y_1-y_2\|,$
$\forall \,\, y_1,\,y_2 \in M.$
\par (iv)\,$d(Vy,y)=\|Vy-y\|=\|2(m-y)\|\geq K\|m-y\|$ with
$K=\alpha_F(r)>1$. \par From (i)-(iv) we have shown the conditions
of Lemma 1 are satisfied, so $m=T(\frac{a-b}{2})$ is a fixed point
for every surjective isometry of $M$. \par Since
$2T(\frac{a-b}{2})\in T(E)$, we can find an $x_0\in E$ such that
$$T(x_0)=2T(\frac{a-b}{2})=2m.\eqno{(2)}$$
\par Define $S: M\rightarrow M\, \mbox {by} \,S(y)=T(x_0-T^{-1}(y))$.
We shall prove that $S$ is well-defined and it is an isometry from M
onto M.
\par (i)\,$S$ is well-defined.\par If $T(x_1)=T(x_2)=y$, then
$$\|T(x_0-x_1)-T(x_0-x_2)\|=\phi\|x_1-x_2\|=\|Tx_1-Tx_2\|=0.$$
\par(ii)\,$S$ is an isometry from M onto M.\par
For any $y_1,y_2\in M$ with $T(x_1)=y_1 \,\mbox {and}\, T(x_2)=y_2$,
we have
$$
\|Sy_1-Sy_2\|=\|T(x_0-x_1)-T(x_0-x_2)\|=\phi(\|x_1-x_2\|)=\|Tx_1-Tx_2\|=\|y_1-y_2\|.$$\par
For any $y\in M\subseteq T(E)\,\mbox {with}\,T(x)=y$, we have
$$ S(y)=T(x_0-x)\, \mbox{and}\,\|Tx\|=\|2m-Tx\|. \eqno {(3)}$$By (2)
and (3), we have
\begin{eqnarray*}
\|Sy\|&=&\|T(x_0-x)\|=\|T(x_0-x)-T(0)\|=\phi(\|x_0-x\|)=\|Tx_0-Tx\|\\
       &=&\|2m-Tx\|=\|Tx\|=\|T(x)-T(0)\|=\phi(\|x-0\|)=\phi(\|x_0-(x_0-x)\|)\\
       &=&\|T(x_0)-T(x_0-x)\|=\|2m-Sy\|.
\end{eqnarray*}
This implies $S(M)\subseteq M$ and by the definition of $S$, we can
easily get $S^2=id_M$, so $S(M)=(M).$ Thus $m$ is a mixed point of
$S$, i.e. $$ T(\frac {a-b}{2})=m=Sm=T(x_0-\frac {a-b}{2}).
\eqno{(4)} $$ By $(2)$ and (4), we have
\begin{eqnarray*}
\|T(a-b)-2T(\frac{a-b}{2})\|&=&\|T(a-b)-Tx_0\|=\phi(\|a-b-x_0\|)\\
                            &=&\phi(\|\frac{a-b}{2}-(x_0-\frac
                            {a-b}{2})\|)\\
                            &=&\|T(\frac{a-b}{2})-T(x_0-\frac {a-b}{2})\|=0.
\end{eqnarray*}
This implies$$T(a-b)=2T(\frac{a-b}{2}).$$ For the fixed $b$, define
$T_b: E\rightarrow F$ by $T_b(x)=T(x+b)-T(b)$ and $T_b$ has the
following properties:\par (i)\,$T_b$ is a continuous map which
preserves equality of distance with $ T_b(0)=0$ and It has the same
gauge function $\phi$ with $T$.
\par
This follows from
\begin{eqnarray*}
\|T_b(x)-T_b(y)\|&=&\|T(x+b)-T(b)-T(y+b)+T(b)\|\\
                 &=&\|T(x+b)-T(y+b)\|=\phi(\|x-y\|).
\end{eqnarray*}
\par
(ii)\,$2T_b(x)-T_b(y)\in T_b(E), \quad \forall\, x,\,y \in E$
\par For any $x,\,y\in E$, we have
$$2T_b(x)-T_b(y)=2T(x+b)-T(y+b)-T(b).$$
Since $2T(x+b)-T(y+b)\in E$, there is a $ z_0\in E$ such that
$T(z_0)=2T(x+b)-T(y+b)$. Namely
$$2T_b(x)-T_b(y)=T(z_0)-T(b)=T(z_0-b+b)-T(b)=T_b(z_0-b)\,\in
T_b(E).$$ \indent (iii)\,$\|T_b(\frac
{a-b}{2})\|=\phi(\|\frac{a-b}{2}\|)=\|T( \frac{a-b}{2})\|\leq \frac
{r}{4}.$\par This follows from (i).
\par By (i), (ii) and (iii), we have showed that $T_b$ has the same
conditions with $T$. Thus , similarly, we get
$$T_b(a-b)=2T_b(\frac{a-b}{2}).\eqno {(5)}$$  By (5), we have
\begin{eqnarray*}
T(\frac{a+b}{2})&=&T(\frac {a-b}{2}+b)-T(b)+T(b)\\
                &=&T_b(\frac {a-b}{2})+T(b)=\frac
                {T_b(a-b)}{2}+T(b)\\
                &=&\frac{T(a-b+b)-T(b)}{2}+T(b)=\frac{T(a)+T(b)}{2}.
                \end{eqnarray*}
Thus we have $$2T(\frac{a+b}{2})=T(a)+T(b),\quad \mbox {\,for all\,
} \, a,\,b \in E\,\mbox{with}\,\|\frac {a-b}{2}\|\leq\delta.\eqno
(6)$$  \indent
 \textbf{Step 2.}\,Consider the case of $\|\frac
{a-b}{2}\|>\delta$. Since $\|.\|$ is continuous for the
multiplication by scalar, there is a positive integer N such that
$$\|\frac {a-b}{2N}\|\leq\delta.\eqno{(7)}$$
Next we shall show by induction that $$
2T(\frac{a+b}{2})=T(\frac{a+b}{2}+\frac
{n(a-b)}{2N})+T(\frac{a+b}{2}-\frac {n(a-b)}{2N}),\qquad \forall \,
n\geq 1.\eqno {(8)} $$ Putting $n=1$ in (8), by (6) and (7), we find
that (8) holds for $n=1$. Let us suppose (8) holds for $n\leq k$,
then we get $$2T(\frac{a+b}{2})=T(\frac{a+b}{2}+\frac
{k(a-b)}{2N})+T(\frac{a+b}{2}-\frac {k(a-b)}{2N}).\eqno{(9)}$$
$$ 2T(\frac{a+b}{2})=T(\frac{a+b}{2}+\frac
{(k-1)(a-b)}{2N})+T(\frac{a+b}{2}-\frac {(k-1)(a-b)}{2N}).\eqno
{(10)}$$ and by (6) and (7), similarly, we can get
$$2T(\frac{a+b}{2}+\frac {k(a-b)}{2N})=T(\frac{a+b}{2}+\frac
{(k+1)(a-b)}{2N})+T(\frac{a+b}{2}+\frac {(k-1)(a-b)}{2N}).\eqno
{(11)}$$$$2T(\frac{a+b}{2}-\frac {k(a-b)}{2N})=T(\frac{a+b}{2}-\frac
{(k+1)(a-b)}{2N})+T(\frac{a+b}{2}-\frac {(k-1)(a-b)}{2N}).\eqno
{(12)}$$Combining (9), (10), (11) and (12), we get
$$2T(\frac{a+b}{2})=T(\frac{a+b}{2}+\frac
{(k+1)(a-b)}{2N})+T(\frac{a+b}{2}-\frac {(k+1)(a-b)}{2N}).$$ Thus we
have proved (8). Putting $n=N$ in (8), we get $$ T(\frac
{a+b}{2})=\frac {T(a)+T(b)}{2}.$$ \par From step 1 and step 2, we
get $$T(a)=T(\frac{2a+0}{2})=\frac {T(2a)+T(0)}{2}=\frac{T(2a)}{2}
\quad \forall \, a\in E,$$ and $$ T(a+b)=T(\frac {2a+2b}{2})=\frac
{T(2a)+T(2b)}{2}=T(a)+T(b),\quad \forall \,a, b\in E.$$ \indent Thus
T is additive. Since it is continuous and the spaces are real, it is
a linear operator and the proof is complete. \qquad $\Box$
\par
\textbf {Remark 1.}\, Theorem 2.2 shows that the assumption in [10]
which is $ \frac {T(x)+T(y)}{2}\in T(E),$ for all $x,y \,\in E$ can
be dropped and the assumptions about the gauge function $\phi$ for
$T$ can be weakened a lot. \par
 \textbf {Remark 2.}\, If $T(E)$ is a
additive group of $F$ and keeps the other hypothesis of Theorem 2.1,
then T is linear.
\par
 \textbf {Corollary 2.3.}\, Let $E,\,F$ be $p$-normed spaces ($0<p\leq 1$),
Let $T: E\rightarrow F$ with $T(0)=0$ be a continuous map which
preserves equality of distance and $$2T(x)-T(y)\in T(E), \quad
\forall \, x,y \in E,$$  then
\par (i)\,$T$ is linear, \par (ii)\,$T=\beta V $ where $\beta$ is a
real number and $V$ is an isometry from $E$ to $F$.
\par
\textbf {proof.} Since (i) is a result of Theorem 2.1, we just need
 to prove (ii). Let $\phi$ be the gauge function for $T$. By the linearity of $T$, we get
$$ \|Tx-Ty\|=\phi(1)\|x-y\|, \quad  \forall \, x,y\in E.\eqno {(13)}$$ \indent If
$\phi(1)=0$, then $T\equiv0$. Let $\beta=0$, then (ii) holds for any
isometry $V$ from $E$ to $F$.\par If $\phi(1)\neq 0$, then let
$\beta=\pm(\phi(1))^\frac {1}{p} \,\mbox{and} \, V=T/\beta.$ By (13)
and the linearity of $T$, we have
\begin{eqnarray*}
\|Vx-Vy\|&=&\|\frac{T}{\beta}(x)-\frac{T}{\beta}(y)\|=\frac{1}{\phi(1)}\|Tx-Ty\|\\
         &=&\frac{1}{\phi(1)}\phi(1)\|x-y\|=\|x-y\|.\quad  \forall
        \,x,\,y \in E.
 \end{eqnarray*}
This implies $V$ is an isometry and the proof is complete. \qquad
$\Box$ \par In particular, when $\phi$ is the identity function, i.e
$\phi(t)=t,$ for all $ t\geq 0$. By Theorem 2.2, we get the
following result which gives a condition for the linearity of not
necessarily surjective isometries between two $F^*$-spaces.
\par \textbf {Corollary 2.4.}\, Let $E,F$ be $F^*$-spaces, and there exists a positive number $r$
such that $ \alpha_F(r)=\inf\{\frac{\|2x\|}{\|x\|}:\, x\in F,\,
0<\|2x\|\leq r\}>1$, Let $V: E\rightarrow F $ with $V(0)=0$ be an
isometry, i.e. $$\|Vx-Vy\|=\|x-y\|\quad \forall \, x,y \in E.$$ and
satisfy $$2V(x)-V(y) \in V(E) \quad \forall \,x,\,y\in E.$$ Then $V$
is linear.\par The following remark was also shown in [10].
\par \textbf {Remark 3.}\, Let $E$, $F$ be normed spaces. Let $V: E\rightarrow F$ be an isometry
with $V(0)=0$. If $F$ is strictly convex, then $$2V(x)-V(y)\in V(E),
\quad \forall x,\,y\in E$$ (Thus $V$ is linear).
\par
\textbf{Proof.}\, For any $x,\,y\in E$, we have
$$\|V(2x-y)-V(x)\|=\|x-y\|=\|Vx-Vy\|,\eqno {(14)}$$
$$\|V(2x-y)-V(y)\|=2\|x-y\|=2\|Vx-Vy\|.\eqno {(15)}$$
By (14) and (15), we get
\begin{eqnarray*}
2\|V(x)-V(y)\|&=&\|V(2x-y)-V(y)\|=\|V(2x-y)-V(x)+V(x)-V(y)\|\\
          &\leq&\|V(2x-y)-V(x)\|+\|V(x)-V(y)\|=2\|V(x)-V(y\|.
 \end{eqnarray*}
Thus $$\|V(2x-y)-V(x)+V(x)-V(y)\|=\|V(2x-y)-V(x)\|+\|V(x)-V(y)\|.$$
Since $F$ is strictly convex, there is a $\lambda >0$ such that
$$ V(2x-y)-V(x)=\lambda \{V(x)-V(y)\}.$$
By (14), we have $\lambda=1.$ Namely $2V(x)-V(y)=V(2x-y)\in
V(E).$\qquad $\Box$ \par By the remark, we can consider the result
of Baker [11] on isometries in strictly convex normed spaces as a
consequence of corollary 2.4.
\par \textbf {Remark 4.}\, Following [13], we give an example to show that
Corollary 2.4 is no longer true if $V$ is an isometry without the
hypothesis $2V(x)-V(y)\in V(E),$ for all $x,\,y \in E.$ \par Let
$E=(R^2,\,\|.\|_p),\,F=(R^3,\,\|.\|_p)$ be $p$-normed spaces
($0<p\leq1$) equipped with $p$-normed defined by
$$\|(\xi,\,\eta)\|_p=(\mbox {max}(|\xi|,\,|\eta|))^p,\quad\forall
\,x=(\xi,\,\eta)\in R^2,$$$$\|(\xi,\,\eta,\,\zeta)\|_p=(\mbox
{max}(|\xi|,\,|\eta|,\,|\zeta|))^p,\quad\forall
\,x=(\xi,\,\eta,\,\zeta)\in R^3.$$ \indent Let $V :E\rightarrow F$
be defined by $V(\xi,\,\eta)=(\sin\xi,\,\xi,\,\eta).$ It is easy to
verity that $V$ is an isometry which is not linear, and if we let
$\xi_1=\frac{\pi}{2},\,\xi_2=\frac{\pi}{3}$, we find that
$2V(\xi_1,\,\eta)-V(\xi_2,\,\eta)\overline{\in}\,V(E),\quad \forall
\,\eta \in R.$\\
\\ \par  \textbf{Acknowledgment}\\  \par The author is grateful to Professor
Ding Guanggui for his encouraging advice and many useful
suggestions.


{\small Tan Dongni
\par
  School of Mathematics Science
  \par
  Nankai University
  \par Tianjin300071, China}\baselineskip 5mm

\end{document}